\documentclass[draft,twoside]{article}
\newcommand{\mytitle}{Equidistribution of Elements of Norm 1 in Cyclic
Extensions}
\newcommand{\keywords}{equidistribution, cyclic number field, Hecke
  zeta function, norm 1}
\newcommand{\msc}{
11K36  %Number Theory Well-distributed sequences and other variations
11R42  %Zeta functions and $L$-functions of number fields
11R04  %Algebraic numbers; rings of algebraic integers
11R27  %Units and factorization
}

\usepackage{amsmath,amssymb,
  amsthm,amsfonts,amscd,mathrsfs,pifont,upgreek,color}
\usepackage{eucal}  % Euler fonts
\usepackage{verbatim}
\usepackage{ifpdf}
\usepackage[colorlinks=true, citecolor=blue]{hyperref}

\ifpdf
        \usepackage[pdftex]{graphicx}
\else
        \usepackage[dvips]{graphicx}
\fi
\newtheorem{thm}{Theorem}[section]

\newtheorem{lemma}[thm]{Lemma}
\newtheorem{prop}[thm]{Proposition}

\newtheorem{thm*}{Theorem}[]
\newtheorem{cor*}[thm*]{Corollary}
\newtheorem{claim*}[thm*]{Claim}
\newtheorem{lemma*}[thm*]{Lemma}
\newtheorem{prop*}[thm*]{Proposition}
\newtheorem{conj*}[thm*]{Conjecture}

\theoremstyle{definition}

\newtheorem{problem*}{Problem}[section]

\newtheorem{question*}{Question}[section]

\newtheorem{defn*}{Definition}

\theoremstyle{remark}

\newcommand{\qq}[1]{\qquad \mbox{#1} \qquad}

\newcommand{\BB}[1]{\ensuremath{\mathbb{#1}}}

\newcommand{\N}{\ensuremath{\BB{N}}}
\newcommand{\Q}{\ensuremath{\BB{Q}}}
\newcommand{\R}{\ensuremath{\BB{R}}}
\newcommand{\Z}{\ensuremath{\BB{Z}}}
\newcommand{\C}{\ensuremath{\BB{C}}}
\newcommand{\T}{\ensuremath{\BB{T}}}

\newcommand{\mf}{\ensuremath{\mathfrak}}

\newcommand{\rp}[1]{\mathrm{Re}(#1)}

\newcommand{\piecewise}[1]{
\left\{
\begin{array}{ll}
#1
\end{array}
\right.
}

\newtheorem*{theorem*}{\rm \bf{Main Theorem}}

\numberwithin{equation}{section}

\numberwithin{equation}{section}

\usepackage{times}

\begin{document}
\title{\bfseries\sffamily \mytitle}
\author{{\sc Kathleen L.~Petersen} and {\sc Christopher D.~Sinclair}}
\maketitle

\begin{abstract}
Upon quotienting by units, the elements of norm 1 in a number field $K$
form a countable subset of a torus of dimension $r_1 + r_2 - 1$ where
$r_1$ and $r_2$ are the numbers of real and pairs of complex
embeddings.  When $K$ is Galois with cyclic Galois group we
demonstrate that this countable set is equidistributed in a finite cover of this torus
with respect to a natural partial ordering induced by Hilbert's Theorem 90.
\end{abstract}

{\bf MSC2010:} \msc

\vspace{1cm}

\section{Introduction}

Our main theorem is the following equidistribution result.

\begin{theorem*}
The elements of norm 1 in a cyclic number field $K$ are contained in a torus under a  natural quotient.  In this quotient, these elements are equidistributed with respect to a so-called visible height function $h$ induced by Hilbert's Theorem 90.
\end{theorem*}
This is stated precisely in Theorem~\ref{main-thm} in Section~\ref{section:heighthilbert} and generalizes previous work on the case of quadratic extensions \cite{MR2854219}.

Let $K$ be a fixed embedding of a number field of degree $d$ over $\Q$.
We define $\mathcal O$ to be the ring of integers of $K$ and $\mathcal U$ the group of units. Let $\mathcal N$ be the subset of elements in $K$ with norm ($N = N_{K/\Q}$) equal to $1$. By Hilbert's Theorem 90, if $K$ is cyclic with Galois group generated by (a preferred generator) $\sigma$, for $\beta \in \mathcal N$, there exists $\alpha \in \mathcal O$ such that $\beta = \alpha/\sigma(\alpha)$.  We define
\[
\pi:
K^{\times} \rightarrow \mathcal N \qq{by} \pi(\alpha) =
\alpha/\sigma(\alpha).
\]
The map $\pi$ is a surjective homomorphism with kernel $\Q^{\times}$, so that the induced map $\tilde{\pi}:K^{\times}/\Q^{\times} \rightarrow \mathcal N$ is an isomorphism.
We will say $\alpha \in \mathcal O$ is a {\em visible point} for $\beta$ if $\pi(\alpha) = \beta$  and $|N(\alpha)|$ is minimal for all integers with this property. In this case $\tilde{\pi}^{-1}(\beta)$ is the coset $\alpha \thinspace \Q^{\times}.$ We denote the set of visible points by $\mathcal V$.  We then define a visible height  on $\mathcal N$ by  $h(\beta) = |N(\alpha)|$ for $\alpha \in \pi^{-1}(\beta)$ a visible point.

The quotient $K^{\times}/\Q^{\times}$ decomposes as $\text{Tor}(K^{\times}/\Q^{\times})\oplus F_{\infty}$ where the first factor is finite and the second is the free group with countably infinite rank. By the isomorphism $\tilde{\pi}$ we conclude that the group structure of $\mathcal N$ is the product of the group of roots of unity in $K$ and a free group.

We consider a natural quotient of $\mathcal N$  into a compact group to formulate our equidistribution question using Weyl's Criterion.  Specifically, with $\log$ denoting the regulator map (see Equation~\ref{regulator}) we consider equidistribution of $\overline{\mathcal N} = \log \mathcal N + \log \pi(\mathcal U)$ where $\mathcal U$ is the group of units in $\mathcal O$.  The set $\pi(\mathcal U)$ is a finite index subgroup of $\mathcal U$  (as proven in Proposition~\ref{torusprop}), and from Dirichlet's Unit Theorem it follows that $\overline{\mathcal{ N}}$ is contained in a torus.  In Proposition~\ref{main-lemma} we show that $h$  induces a well-defined height, which we also call $h$, on $\overline{\mathcal N}$.

To prove that $\overline{\mathcal{ N}}$  is equidistributed with respect to $h$, we use Weyl's criterion and the Wiener-Ikehara Tauberian Theorem to reduce the problem to an analysis of $L$-functions.  After corresponding our natural $L$-series to a partial Hecke $L$-series, we  use known results about the analytic properties of partial Hecke $L$-series to demonstrate that Weyl's Criterion is satisfied.

\subsection{Equidistribution and Weyl's Criterion}

Given a probability space $(X, \mathcal F, \nu)$ and countable set $C \subset X$ we call $h: C \rightarrow [0, \infty)$ a {\em height} on $X$ if
\[
C(r) := h^{-1}[0,r] = \{ x \in C : h(x) \leq r \}
\]
 is finite for every $r > 0$. We then say that $C$ is {\em equidistributed} (or uniformly distributed) in $(X, \mathcal F, \nu)$ with respect to the height $h$ if for any measurable set $B \in \mathcal F$ with positive measure,
\[
\lim_{r \rightarrow \infty} \frac{\# B \cap C(r)}{\# C(r)} = \nu(B).
\]
In words, $C$ is equidistributed with respect to $h$ if, for every positive $B \in \mathcal F$ the proportion of points of bounded height which lie in $B$ converges to the measure of $B$ as the height bound increases.

When $X$ is a compact abelian group with Borel $\sigma$-algebra and Haar probability measure---we say simply that $C$ is equidistributed in $X$ with respect to $h$. In this case, we may use Weyl's criterion to establish equidistribution. Specifically, $C$ is equidistributed in a compact abelian group $X$ with respect to $h$ if, for any character $\chi: X \rightarrow \T$,
\[
\lim_{n \rightarrow \infty} \frac{1}{\# C(r)} \sum_{x \in
  C(r)} \chi(x) = \piecewise{1 & \mbox{if $\chi$ is
    trivial;} \\ 0 & \mbox{otherwise}.}
\]

We now precisely define the compact abelian group and height to which we will apply Weyl's Criterion.

\subsection{A Hilbert-Dirichlet Torus}

For a number field $K$, the set of archimedean places is denoted $S_{\infty}$, and has
cardinality $r_1 + r_2$ where $r_1$ is the number of real places and
$r_2$ is the number of complex places.  For each place $v \in
S_{\infty}$ we denote by $\| \cdot \|_v$ either the usual absolute
value if $v$ is real, or the usual complex absolute value squared if $v$ is
complex.

We define the {\em regulator map} $\log: K^{\times} \rightarrow \R^{r_1 + r_2}$ by
\begin{equation}
\label{regulator}
\alpha \mapsto \left( \log \| \alpha \|_v \right)_{v \in S_{\infty}}.
\end{equation}
We also define $\Sigma: \R^{r_1 + r_2} \rightarrow \R$ by $\mathbf x
\mapsto x_1 + x_2 + \cdots + x_{r_1 + r_2}$.  Then, by the proof of
Dirchlet's Unit Theorem, $\log \mathcal U$ is a lattice in $\ker \Sigma$ and
the kernel of $\log$ restricted to $\mathcal U$ is the group of roots of unity in $K$.  This is  usually stated as
\[
\mathcal U \cong \mathrm{Tor}(\mathcal U) \times \Z^{r_1 + r_2 - 1}.
\]
The quotient
\[
T = \ker \Sigma/\log \mathcal U
\]
is isomorphic to the torus $\T^{r_1 + r_2 - 1}$.
Our normalization of $\|\cdot \|_v$ together with the fact that $\mathcal N$ consists of elements of norm 1 implies that $\log \mathcal N$  lies in $\ker \Sigma$.
From here forward we will assume that $K$ is Galois over $\Q$ with
cyclic Galois group $G$ generated by $\sigma$.\footnote{Note that this implies
that $K$ is either totally real, or totally imaginary, but we will
still write, for instance, $\mathcal U \cong \mathrm{Tor}(\mathcal U) \times \Z^{r_1 + r_2 - 1}$ with the understanding that one of $r_1$ and $r_2$ is 0.}

\begin{prop}
\label{torusprop}
The set $S := \ker \Sigma / \log \pi(\mathcal U)$ is a torus and is a  finite cover of $T=\ker \Sigma/\log \mathcal U$.
\end{prop}
\begin{proof}
    This will follow from the fact that $\pi(\mathcal U)$ is a finite index subgroup of $\mathcal U$, which we will now establish.  We will show that for all $u \in \mathcal U \cap \mathcal N$, $u^e\in \pi(\mathcal U)$ where
    $e$ is the least common multiple of the ramification indices of all  prime ideals in $\mathcal O$.

Since $u\in \mathcal N$,  there is a visible point $\alpha$ such that $u = \pi(\alpha)$.
    For each prime ideal $\mf p$, there is a non-negative integer $n_{\mf p}(\alpha)$ such that
    \[
    \alpha \mathcal O = \prod_{\mf p} \mf p^{n_{\mf p}(\alpha)} = \sigma(\alpha) \mathcal O,
    \]
    where the second equation follows because $\alpha/\sigma(\alpha) = u$ is a unit.  If $\alpha$ is a unit then every $n_{\mf p}$ is zero.  If $\alpha$ is not a unit, at least one $n_{\mf p} > 0$ and, in any case, $n_{\mf p} = 0$ for all but finitely many primes $\mf p$.

  The Galois group acts transitively on the primes above a given rational prime $p$. This implies that if $\mf p, \mf q$ are primes above $p$, then $n_{\mf p}(\alpha) = n_{\mf q}(\alpha)$ and in this case we can write $n_p(\alpha)$ for this common integer. That is,
    \[
    \alpha \mathcal O = \prod_p \prod_{\mf p | p} \mf p^{n_p(\alpha)}.
    \]

    We come to the central observation: if $n_p(\alpha) \geq e_p$, the ramification index of $p$,  then $p \mathcal O$ divides $\alpha \mathcal O$ and hence $\alpha = p \beta$ for some $\beta \in \mathcal O$. In this situation $\pi(\alpha) = \pi(\beta)$ and $p^d |N(\beta)| = |N(\alpha)|$ contradicting the visibility of $\alpha$. In particular $n_p(\alpha) = 0$ for all non-ramified primes $p$.

    Consider   $u^e= \pi(\alpha^e)$.  For all primes $p$ we write $e=e_p g_p$,  for some positive integers $g_p$.  We have that
    \[
    \alpha^e \mathcal O = \prod_p \prod_{\mf p | p} \mf p^{en_p(\alpha)} = \prod_p \prod_{\mf p | p} \mf p^{e_p g_p n_p(\alpha)} = \prod_p  p^{ g_p n_p(\alpha)}.
    \]
    As such, $\alpha^e \mathcal O$ is generated by a rational integer.  It follows that  $\alpha^e = n w$ for some $n \in \mathbb N$ and $w\in \mathcal U$.
    We have $u^e =\pi(\alpha^e) = \pi(nw) = \pi(w)$ as needed.
\end{proof}

\subsection{A Visible Height From Hilbert's Theorem 90}\label{section:heighthilbert}

We now prove a key lemma from which it follows that the visible height is well-defined on $\mathcal N$.
\begin{lemma}
\label{claim:2}
Assume that  $\alpha\in \mathcal V$ is a visible point for $\beta \in \mathcal N$, and $\gamma \in \mathcal O$.  Then  $\beta = \pi(\gamma)$ if and only if there exists a non-zero rational integer $n$ such that $\gamma = n \alpha$.  Such an $n$ and $\alpha$ are unique up to multiplication by $-1$.
\end{lemma}
In particular, $\alpha$ is a visible point for $\beta$ if and only if $-\alpha$ is.
\begin{proof}
    If $\gamma = n \alpha$ and $\alpha$ is a visible point for
    $\beta$ we have
    \[
    \beta = \pi(\alpha) = \frac{\alpha}{\sigma(\alpha)} = \frac{n
      \alpha}{ n \sigma(\alpha)} = \frac{n\alpha}{ \sigma(n \alpha) }=
    \pi(n \alpha) = \pi(\gamma).  \]

    Conversely, assume that  $\pi(\alpha) = \pi(\gamma)$ then
    \[
    \frac{\alpha}{\gamma} = \frac{\sigma(\alpha)}{\sigma(\gamma)} =
    \sigma\left( \frac{\alpha}{\gamma} \right),
    \]
    and hence $\alpha/\gamma$ is fixed by $G$.  It follows that there
    exist relatively prime rational integers $m$ and $n$ such that
    $\alpha/\gamma = m/n$, and by the minimality of $|N(\alpha)|$, $|m| \leq
    |n|$.  Thus both $\alpha$ and $\gamma = n \alpha/m$ are algebraic
    integers which map onto $\beta$.

    If we assume $n/m > 0$, then $n/m\geq 1$ since $|m|\leq |n|$ and
    for $j = \lfloor n/m \rfloor$ we have that
    $\gamma - j \alpha$ is an algebraic integer.  If $n/m \neq j$ then it is non-zero with
    \[
    \pi(\gamma - j \alpha) = \pi\left(\left( \frac{n}{m} - j \right)
    \alpha \right)  = \frac{(n/m - j) \alpha}{\sigma((n/m -
    j)\alpha)} = \pi(\alpha).
    \]
Since $n/m-j \in \mathbb Z - \{0\}$,
     \[
    |N(\gamma - j \alpha)| = \left|N\left(\Big(\frac{n}{m} - j\Big)
        \alpha\right)\right| = \left(\frac{n}{m} - j\right)^d |N(\alpha)|
    < |N(\alpha)|.
    \]
    This contradicts the minimality of $|N(\alpha)|$.  Therefore,  $n/m= j$ so that  $n/m$ is an integer and
     $\gamma = n \alpha$.
    The case when $n/m < 0$ is similar.

    To show uniqueness, by the above if both $\alpha$ and $\alpha'$ are visible points for $\beta$ then $\pi(\alpha)=\pi(\alpha')$ and  we have $\alpha=n' \alpha'$ and also $\alpha' = n \alpha$ for $n, n'\in \mathbb{Z}$. But then $n, n'=\pm1$.
\end{proof}

Hilbert's Theorem 90 implies that for every $\beta \in \mathcal N$
there is a visible point $\alpha \in \mathcal V$, and in this situation, we define
\[
h(\beta) = |N(\alpha)|.
\]
For any $n > 0$ there are finitely many integers in $\mathcal O$ with absolute norm bounded by $n$. It follows from this  fact and Lemma~\ref{claim:2}  that $h$ is a well-defined height on $\mathcal N$. A height on $\mathcal N$ does not necessarily induce a well defined on $\overline{\mathcal N}$. We observe that for each $\beta \in \mathcal N$ and $u \in \pi(\mathcal  U)$, $h(u \beta) = h(\beta)$.  By the uniqueness statement in Lemma~\ref{claim:2} we arrive at the following.

\begin{prop}
\label{main-lemma}
The height $h$ induces a well defined height on $\overline{\mathcal N} := \log \mathcal N + \log \pi(\mathcal  U)$.
\end{prop}

The set $\mathcal N$ of norm 1 elements is contained in $\ker \Sigma$ and so $\overline{\mathcal N} = \log \mathcal N + \log \pi(\mathcal  U)$ is contained in $\ker \Sigma/\log \pi(\mathcal U).$
Propositions~\ref{torusprop} and \ref{main-lemma} imply that $\overline{\mathcal N}$ is a countable set in a compact abelian group with a natural height, and therefore we can state our Main Theorem precisely.

\begin{thm}
\label{main-thm}
$\overline{\mathcal N}=\log \mathcal N + \log \pi(\mathcal U)$ is equidistributed in $S=\ker \Sigma/\log \pi(\mathcal U)$ with respect to $h$.
\end{thm}

\subsection{Hecke $L$-Series and the Wiener-Ikehara Tauberian Theorem}
Let $\overline{\mathcal{N}}(r)=\{x\in \overline{\mathcal N} : h(x) \leq r\}.$ The function
\[
r \mapsto \sum_{x \in \overline{\mathcal{N}}(r)} \chi(x)
\]
is the summatory function of the $L$-series
\[
L(\chi; s) = \sum_{x \in \overline{\mathcal{N}}} \frac{\chi(x)}{h(x)^s}.
\]
This observation is useful since asymptotics of the summatory
function (as a function of $r$) follow from analytic properties of $L(\chi; s)$ using the Wiener-Ikehara Tauberian theorems.  Specifically, if
we can show that, for the trivial character $\chi = 1$, $L(s; 1)$
has a pole at $s = \sigma_0$ and there exists $\epsilon > 0$ such that
for all other characters, $L(\chi; s)$ is analytic for $\rp{s} >
\sigma_0 - \epsilon$, then there exists a non-zero constant $C$ such
that as $r \rightarrow \infty$,
\[
\# \overline{\mathcal{N}}(r) \sim C r^{\sigma_0} \qq{and} \sum_{x \in \mathcal{\overline N}(r)}
\chi(x) = o(r^{\sigma_0}),
\]
and Weyl's criterion will be satisfied.  See \cite[Ch.VII \S3]{MR0160763} or \cite[Ch. XV \S3,5]{MR1282723}.

Let $\mf O$ be the set of non-zero principal integral ideals in $K$ and let $\N$ denote the ideal norm. Let $\iota: K^{\times} \rightarrow (\R^\times)^{r_1} \times (\C^\times)^{r_2}$ be a preferred Minkowski embedding of $K^{\times}$.  Given a character, we denote the complex conjugate character with a bar.

\begin{thm}
  \label{hecke-form}
Each character $\chi$ on $S$ lifts to a continuous character $\chi_{\infty}$ on $(\R^\times)^{r_1} \times (\C^\times)^{r_2}$. Define $\chi_{\mathfrak I}$ to be the character on the the group of principal fractional ideals $\mathfrak I$ given by
\[
\chi_{\mathfrak I}(\alpha \mathcal O) = \overline{\chi}_{\infty}(\iota(\alpha)).
\]
Then, in the half-plane $\rp{s}>1$ the function $L(\chi;s)$ converges absolutely and
\[
L(\chi; s) = \frac{1}{\zeta(ds)} \sum_{\mf a \in \mf O} \frac{\overline{\chi}_{\mathfrak I}(\mf a)}{\N \mf a^s},
\]
where  $\zeta(s)$ is the Riemann zeta function.  Moreover, $\overline{\chi}_{\mathfrak I}$ is trivial if and only if $\chi$ is trivial. For $\chi$ non-trivial $L(\chi;s)$ has an analytic continuation to the half-plane $\rp{s}>\tfrac1d$, and $L(\chi;1)$ has a meromorphic continuation to this half-plane with a simple pole at $s=1$.
\end{thm}
Our notation is suggestive as $\chi_{\infty}$ extends to the $\infty$-type of the unramified Hecke character $\chi_{\infty} \overline{\chi}_{\mathfrak I}$ on the id\`ele class group  of $K$. It will then follow that
\[
\zeta_K(\chi_{\mathfrak I}; s) := \sum_{\mf a \in \mf O} \frac{\overline{\chi}_{\mathfrak I}(\mf a)}{\N \mf a^s},
\]
is a {\em partial} Hecke $L$-series (here {\em partial} means summed over the class of principal ideals).  If $K$ happens to have class number 1, then $\zeta_K(\chi_{\mathfrak I}; s)$ is a (complete) Hecke $L$-series.

In \cite{MR1544392} it is proved that
$\zeta_K(\chi_{\mathfrak I}; s)$ has an analytic continuation (as a function of
$s$) to the entire plane, except when $\chi_{\mathfrak I}=1$ is trivial, in which
case it has a meromorphic continuation to a function with a single
simple pole at $s=1$.  It follows that, if $\chi$ is nontrivial, then
so is $\chi_{\mathfrak I}$, and $L(\chi; s) = \zeta_K( \chi_{\mathfrak I}; s)/ \zeta(d s)$
has an analytic continuation to the half-plane $\rp{s} > 1/d$.  If
$\chi=1$ is trivial, then $L(\chi; s)$ has a meromorphic continuation to
the same half-plane, but with a simple pole at $s = 1$.  Weyl's criteria now follows from the Wiener-Ikehara Tauberian Theorem.  Therefore, Theorem~\ref{main-thm} follows from Theorem~\ref{hecke-form}.

% \subsection{The Ideal Height}

% Before embarking on the proofs, we remark that there is a natural height on $\wt{\mathcal N}$ that can be defined for arbitrary numbers fields. If $\mf a$ is a fractional ideal of $\mathcal O$, then we may factor $\mf a$ into prime ideals
% \[
% \mf a = \prod_{\mf p} \mf p^{v_{\mf p}}
% \]
% where the product is over all prime ideals of $\mathcal O$ and $v_{\mf p}$ is an integer which is non-zero for ony finitely many primes $\mf p$. The ideal height is then defined by
% \[
% j(\mf a) = \prod_{\mf p} \N \mf p^{|v_{\mf p}|}.
% \]
% We may then define the height $j: \wt{\mathcal N} \rightarrow \N$ by $j(\wt \xi) := j(\xi \mathcal O)$. This is well defined because if $\zeta$ is another generator of $\xi \mathcal O$, then $\zeta = u \xi$ for some unit $u \in \mathcal U$ and hence $\wt \zeta = \wt \xi$. We will consider equidistribution of $\wt{\mathcal N}$ with respect to this height $j$ in a subsequent paper.

\section{The Proof of Theorem~\ref{hecke-form}}

Starting with a continuous character $\chi$ on $S$ we note that $\overline{\mathcal N}=\log \mathcal N + \log \pi(\mathcal U)$ is dense in $S$.   Given $\beta  \in \mathcal N$,
\[
\pi^{-1}(\beta) = \{ x \in K^{\times} : x/\sigma(x) = \beta \}.
\]
and $K^{\times}$ is the disjoint union
\[
K^{\times} = \bigsqcup_{\beta \in \mathcal N} \pi^{-1}(\beta).
\]
We then define $\chi_K : K^{\times} \rightarrow \T$ by $\chi_K(x) = \chi(\log \beta + \log \pi (\mathcal U))$ whenever $x \in \pi^{-1}(\beta)$. 
With Lemma~\ref{claim:2} then  $\chi_K(r)=1$ for all $r\in \Q$, which corresponds to the fact that if $\pi(\alpha) = \beta$ then  $\tilde{\pi}^{-1}(\beta) = \alpha \Q^{\times}$.  
Note that if $u$ is a unit, $\chi_K(u) = \chi(\log \pi(\mathcal U)) = 1$.

 For all $x\in K^{\times}$ we define   $\chi_{\infty}(\iota(x)) = \chi_K(x)$.
Because $\iota(K^{\times})$ is dense in $(\R^{\times})^{r_1} \times (\C^{\times})^{r_2}$, the continuity of $\log$ and $\pi$ mean that we may uniquely define
\[
\chi_{\infty} : (\R^{\times})^{r_1} \times (\C^{\times})^{r_2} \rightarrow \T
\]
by considering points in $(\R^{\times})^{r_1} \times (\C^{\times})^{r_2}$ as limits of elements of $K^{\times}$ under the (multiplicative) Minkowski embedding $K^{\times} \hookrightarrow (\R^{\times})^{r_1} \times (\C^{\times})^{r_2}$.

It remains to show that $\chi_{\infty}$ is a character. This will follow again by density and the fact that $\chi_K$ is a multiplicative homomorphism. To see this, suppose $x, y \in K^{\times}$. Since $\pi$ is multiplicative $ \pi(x y) =  \pi(x) \pi(y)$, and $\chi$ is a character on $S$ we have that
\begin{align*}
\chi_K(x y) &= \chi(\log \pi(x y) + \log \pi(\mathcal U)) \\
&= \chi(\log \pi(x) + \log \pi(y) + \log \pi(\mathcal U)) \\
&= \chi(\log \pi(x) + \log \pi(\mathcal U)) \chi(\log \pi(y) + \log \pi(\mathcal U)) \\
&= \chi_K(x) \chi_K(y).
\end{align*}
Clearly, $\chi_K$ is trivial if and only if our original character $\chi$ is trivial.

We wish to relate $L(\chi; s)$ to a Hecke $L$-function in order to establish the necessary analytic properties to verify Weyl's criterion. Let $\mf V$ be the set of principal ideals generated by visible points.  If $\mf v \in \mf V$ then there is a visible point $\alpha$ such that $\mf v = \alpha \mathcal O$. If $\alpha'$ is another visible point such that $\mf v = \alpha' \mathcal O$ then $\alpha' = u \alpha$ for some $u \in \mathcal U$. Then,
\[
\log(\pi(\alpha')) = \log(\pi(\alpha)) + \log \pi(u), \qq{and hence} \overline{\pi(\alpha)} = \overline{\pi(\alpha')}
\]
in $\overline{\mathcal N}$
and hence $\log \circ \pi$ induces a bijection from  $\mf V$ to $\overline{\mathcal N}$. We write $\overline{\mf v}$ for the image of $\mf v$ under this map.
This induces a map from $\mathfrak I \rightarrow \overline{\mathcal N}$ by defining for $\mathfrak a = \gamma \mathcal O$ with  $\gamma \in K^{\times}$ that  $\overline{\mathfrak a} = \overline{\pi (\gamma)}$.  As such,
\[
\overline{\chi}_{\mathfrak I}(\mathfrak a) =\chi( \overline{\mathfrak a}) = \chi( \overline{\pi(\gamma)}).
\]
For $\mathfrak v \in \mathfrak V$, as above $h(\overline{\mathfrak v})= |N(\alpha)| = \mathbb N \mathfrak v$  and we have
\[
L(\chi; s) = \sum_{x \in \overline{\mathcal{N}}} \frac{\chi(x)}{h(x)^s} = \sum_{\mf v \in \mf V} \frac{\overline{\chi}_{\mathfrak I}(\mf v)}{\N \mf v^s}.
\]
Multiplying and dividing by the scaled Riemann zeta function $\zeta(d s)$,  up to reordering,
\[
L(\chi; s) = \frac{1}{\zeta(d s)} \sum_{n=1}^{\infty} \sum_{\mf v \in \mf V} \frac{\overline{\chi}_{\mathfrak I}(\mf v)}{n^{d s} \N \mf v^s}
= \frac{1}{\zeta(d s)} \sum_{n=1}^{\infty} \sum_{\mf v \in \mf V} \frac{\overline{\chi}_{\mathfrak I}(\mf v)}{\N (n \mf v)^s}.
\]
By Lemma~\ref{claim:2}, each non-zero algebraic integer $\gamma$ can be written uniquely as $n
\alpha$ for visible point $\alpha$ and positive integer $n$.
The principal ideal    $\mf a = \gamma \mathcal O$ can be written uniquely as $n \mf v$ where $\mf v = \alpha \mathcal O$ and we have $\pi(\gamma) = \pi(\alpha)$ and therefore $\overline{\mf a} = \overline{\mf v}$ and $\overline{\chi}_{\mf I}(\mf a) = \overline{\chi}_{\mf I}(\mf v)$.
 That is, up to reordering,
\[
L(\chi; s) = \frac{1}{\zeta(d s)} \sum_{\mf a \in \mf O}
\frac{\overline{\chi}_{\mathfrak I}(\mf a)}{\N \mf a^s} =\frac1{\zeta(ds)} \zeta_K(\chi_{\mathfrak I}; s).
\]
 The right-hand side of this equation can be written as an absolutely convergent series in the half-plane $\rp{s}>1$ because of the analytic properties of both $\zeta(s)$ and $\zeta_K(\chi_{\mathfrak I}; s)$.  Because this series can be reordered as $L(\chi;s)$ it follows that $L(\chi;s)$  is also an absolutely convergent series in the half-plane $\rp{s}>1$.  The stated analytic properties follow from the analytic properties of $\zeta(s)$ and  $\zeta_K(\chi_{\mathfrak I}; s)$.

\section{Acknowledgments}

We would like to thank Mathematisches Forschungsinstitut Oberwolfach
(Research in Pairs) and its sponsoring agencies for providing a
stimulating and congenial atmosphere in which to work on this
material.  The first author would also like to acknowledge the Simons
Foundation (grant number 430077) for providing partial funding for
this collaboration.

\bibliography{bibliography}

\noindent\rule{4cm}{.5pt}
\vspace{.25cm}

\noindent {\sc \small Kathleen L.~Petersen}\\
{\small Mathematics \& Statistics Department,  University of Minnesota Duluth, Duluth MN 55812} \\
email: {\tt kpete@umn.edu}

\vspace{.25cm}

\noindent {\sc \small Christopher D.~Sinclair}\\
{\small Department of Mathematics, University of Oregon, Eugene OR 97403} \\
email: {\tt csinclai@uoregon.edu}

\end{document}